\documentclass[12pt]{amsart}

\usepackage{graphicx}
\usepackage{caption}
\usepackage{subcaption}
\usepackage{color}
\usepackage{amsmath,amsfonts}
\usepackage{algorithm}
\usepackage{algorithmic}

\theoremstyle{definition}

\theoremstyle{remark}

\begin{document}



\title[Fixed mesh ALE]{Numerical simulation of parabolic\\
      moving and growing  interface problems \\ 
      using small mesh deformation}

\author{Ulrich Langer}
\address{Johann Radon Institute for Computational and Applied Mathematics (RICAM), Austrian Academy of Sciences and Institute of Computational Mathematics, Johannes Kepler University, Altenberger Strasse 69, A-4040 Linz, Austria}
\email{ulanger@numa.uni-linz.ac.at}
\urladdr{http://www.numa.uni-linz.ac.at/~ulanger/}


\author{Huidong Yang}
\address{Johann Radon Institute for Computational and Applied Mathematics (RICAM), Austrian Academy of Sciences, Altenberger Strasse 69, A-4040 Linz, Austria}
\email{huidong.yang@oeaw.ac.at}
\urladdr{http://people.ricam.oeaw.ac.at/h.yang/}



\begin{abstract}
  In this work, we develop a cutting method for solving 
  problems with 
  moving and growing interfaces in 3D. This new method is able to resolve large 
  displacement or deformation of 
  immersed objects
  by combining the Arbitrary Lagrangian-Eulerian 
  method with only small local mesh deformation defined on the 
  reference domain, that is decomposed into the macro-elements. 
  The linear system of algebraic equations arising after the temporal and spatial discretizations 
  of a model parabolic interface heat-conduction-like problem with vector-valued functions 
  is solved by either an all-at-once or a segregated algebraic multigrid method.
\end{abstract}


\maketitle








\section{Introduction}
The conventional Arbitrary Lagrangian-Eulerian (ALE) method (see, e.g., 
\cite{TH:81,JD04:00}) works well for small deformation in many applications. For large 
deformation problem, the ALE method may fail due to the deteriorated mesh quality. 
Some improved ALE methods have been studied, e.g., a method based on the biharmonic extension in 
\cite{Wick11:00}. More work based on 
the so called fixed-mesh ALE approach has been studied, 
e.g., in \cite{NME:NME2740,Codina20091591}. 
The parametric finite element method \cite{SFTR14}, 
the immersed-interface finite element method \cite{YGBLZL08} 
and 
the immersed boundary method \cite{BoffiD11} 
may also be applied in this context. 
An enhanced ALE method combined with the fixed-grid and extended finite element 
method (XFEM) was studied in \cite{FLD:FLD1782}. 
Another promising approach is to use the space-time method, 
that is more flexible to handle 
moving interface problems; see, e.g., \cite{MN11,EKMN15}. 

In this work, we propose an interface capturing method by pre-computing the 
intersection of the moving object immersed in the 
underlying reference tetrahedral elements in three dimension (3D). Combined 
with the ALE method on such reference elements, we are able to deal 
with the moving or growing interface problems with 
large displacement or deformation. In a similar manner as already 
investigated in the earlier work \cite{RBM11,YW11:00,YH11:00}, 
the piece-wise linear finite element basis functions are constructed 
on each macro-element \cite{RBM11}, that is decomposed into four 
pure tetrahedral elements and one octahedral element. 
In addition, the method offers a nice opportunity to keep capturing 
the interface without introducing extra degrees of freedom. 
To test the robustness of the method, we consider a 
model heat-conduction-like problem with vector-valued functions. Such a model 
can be used to handle the mesh movement in the fluid-structure 
interaction simulation, see, e.g., \cite{Razzaq20121156}. 
The construction of robust solution methods for solving 
the arising finite element equations requires 
additional effort. For this, we use both the all-at-once and the segregated methods, that 
employ 
an
algebraic multigrid (AMG) method \cite{FK98:00,UH:02}. 

The remainder of the paper is organized as follows: In Section \ref{sec:model}, we 
set up the model parabolic interface problem. Section \ref{sec:dis} deals with 
the temporal and spatial discretization of the model interface problem. In Section 
\ref{sec:linsol}, we discuss the 
all-at-once and the segregated methods 
for solving the linear system of equations arising from the temporal and spatial discretization. 
We present numerical results of two proposed interface moving problems 
in Section \ref{sec:numa}. Finally, some 
conclusions are drawn in Section \ref{sec:con}. 

\section{A model interface problem}\label{sec:model}
\subsection{Geometrical configurations}
\begin{figure}[htbp]
  \centering
  \scalebox{0.5}{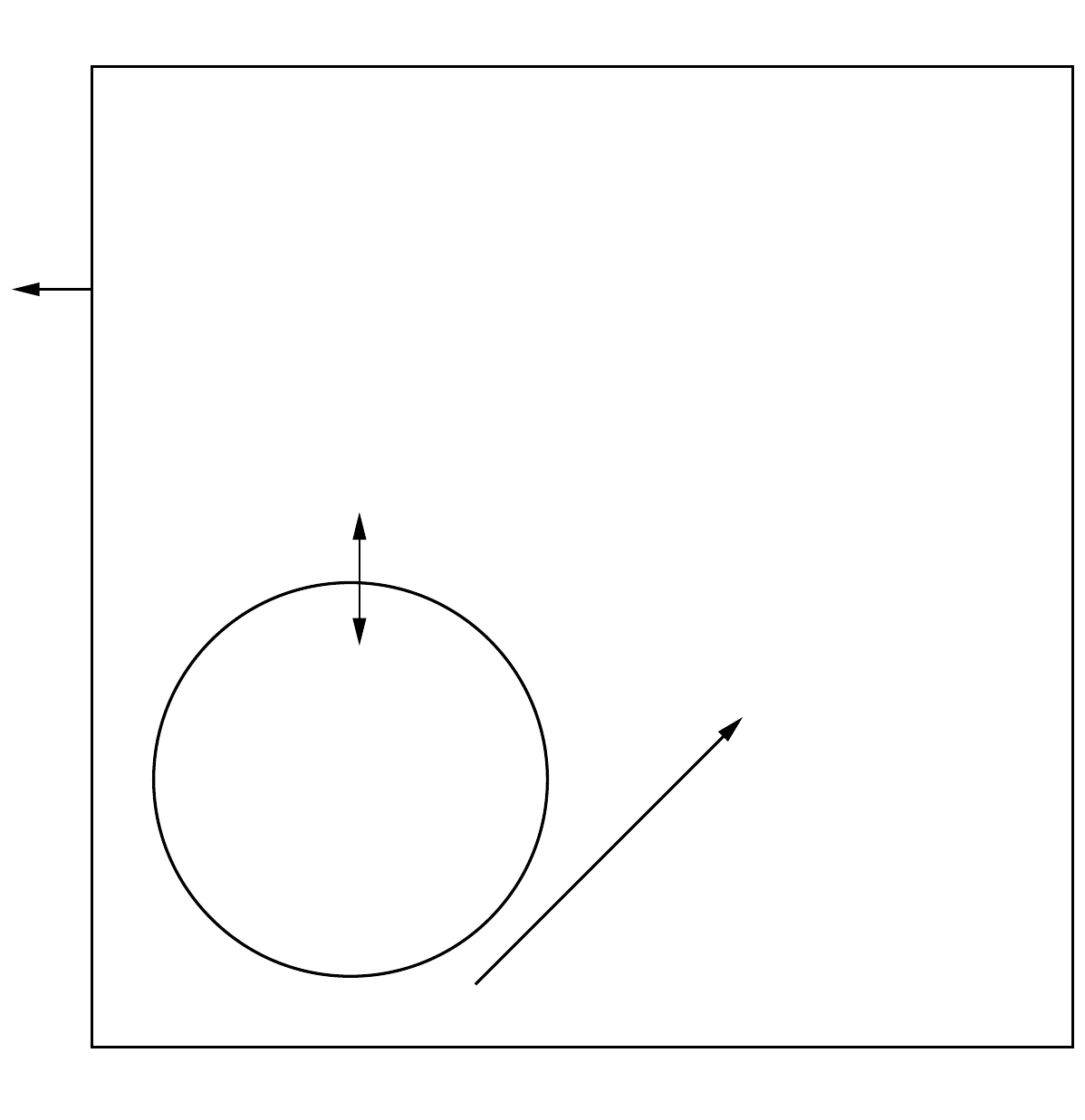}\quad
  \scalebox{0.5}{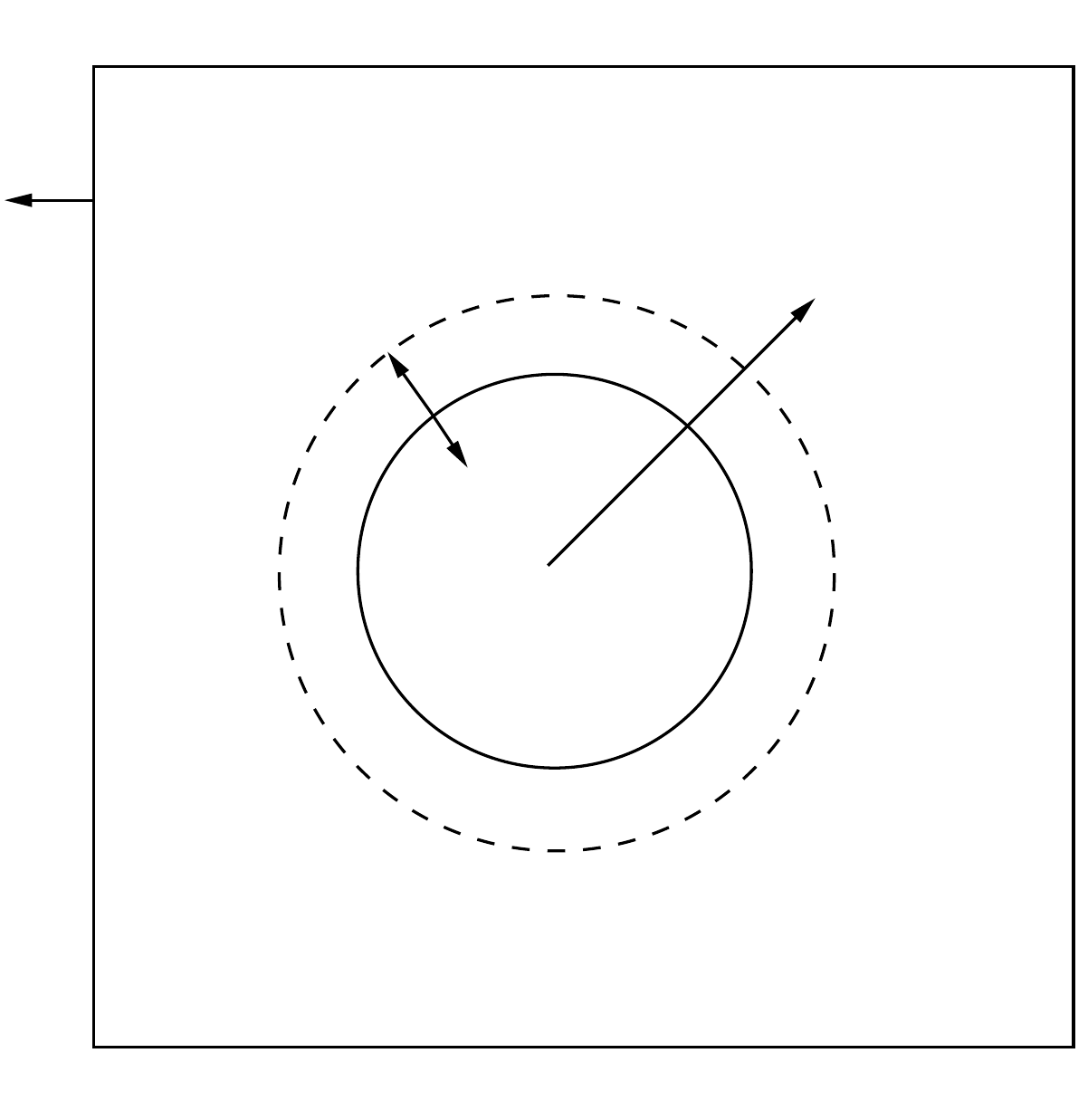}
  \caption{An illustration of two sub-domains for the interface problem: Rigid body 
  motion (left) and growing sub-domain (right).}
  \label{fig:modeldomain}
\end{figure}

We consider a simply connected, 
bounded, polyhedral Lipschitz
domain 
$\Omega\subset {\mathbb R}^3$,
which includes an immersed time-dependent,
sufficiently smooth 
sub-domain $\Omega_1^t=\Omega_1(t)\subset\Omega$, where $t\in I$ denotes the time 
with $I=(0, T]$ being the time interval. The remaining 
sub-domain is $\Omega_2^t=\Omega \backslash \bar{\Omega}_1^t$. 
By $\Gamma^t:=\partial\Omega_1^t\cap\partial\Omega_2^t$, we denote the interface. 
The boundaries of $\Omega$ are denoted by $\Gamma_D$ and 
$\Gamma_N$ such that 
$\partial \Omega=\bar{\Gamma}_D\cup\bar{\Gamma}_N$ and 
$\Gamma_D\cap\Gamma_N=\emptyset$, where 
proper Dirichlet and Neumann boundary conditions are prescribed, respectively. 
We use $n$ to denote the outward unit normal vector on the boundary $\partial\Omega$, 
$n_1$ and $n_2$, the outward unit normal vectors on $\Gamma^t$ with respect to $\Omega_1$ and 
$\Omega_2$, respectively. 
We refer to Fig.~\ref{fig:modeldomain} for an illustration of a sub-domain immersed 
in the big domain. 
We consider two interface problems in this work. 
In the first problem, the sub-domain $\Omega_1$ keeps the shape 
and moves with a constant velocity $v\in {\mathbb R}^3$, i.e., a rigid body motion; 
see the left plot in Fig.~\ref{fig:modeldomain}. 
In the second problem, the sub-domain grows with 
a constant velocity $v$ along the line connecting 
the mass center $p_c$ of the sub-domain and any point $p_b$ 
on the boundary $\partial\Omega_1$; see the right plot in Fig.~\ref{fig:modeldomain}.

\subsection{The model problem with a fixed interface $\Gamma$}
We start to formulate the problem in the fixed sub-domains $\Omega_1$ and 
$\Omega_2$ with proper interface conditions on the fixed interface 
$\Gamma=\partial\Omega_1\cap\partial\Omega_2$. We aim to 
find the solution $u : \Omega \mapsto {\mathbb R}^3$, for all 
$t\in I $, such that
\begin{equation}\label{eq:heatfixedinterface}
  \begin{aligned}
    \partial_t u - \nabla\cdot (a \nabla u ) = 0 & \textup{ in } \Omega_{1}\cup \Omega_{2},\\
            u_1 = u_2  & \textup{ on } \Gamma, \\
            a_1\frac{\partial u_1}{\partial n_1} + a_2\frac{\partial u_2}{\partial n_2}   
            & \textup{ on } \Gamma \\
  \end{aligned}
\end{equation}
with the initial condition $u=0$ at $t=0$, and 
the boundary conditions $u_2 = g_D$ on $\Gamma_D$ and 
$a_2\frac{\partial u_2}{\partial n_2}  = g_N$ on $\Gamma_N$ at $t>0$. 
Here $a=a_1\in {\mathbb R}^{+}$ in $\Omega_1$, 
$a=a_2\in {\mathbb R}^{+}$ in $\Omega_2$, $a_1\neq a_2$, are two different material 
coefficients.  The analysis of such an interface problem with the scalar-valued function 
has been studied, e.g., in \cite{ZJ98,LiZhiLin03,JB96}. In this work, we consider 
the model problem with the vector-valued function, that can be used to model the mesh 
movement in the fluid-structure interaction simulation in our future work. 

\subsection{The model problem with a unfixed interface $\Gamma^t$}
For the interface problem with unfixed interface $\Gamma^t$, 
the time derivative $\partial_t u$ in (\ref{eq:heatfixedinterface}) is not 
well-defined since the computational domain is moving. One of the classical 
approaches is to use the ALE method \cite{TH:81,JD04:00}, in which we 
introduce a displacement defined on the  reference domain $\Omega_R$: 
\[d(x, t):\Omega_R  \times I \mapsto {\mathbb R}^3\] for 
all $x\in\Omega_R$ and $t\in I$, that tracks the monition of the computational 
domain $\Omega$. The ALE mapping ${\mathcal A}^t:\Omega_R\mapsto \Omega^t$ 
for all $t\in I$, where $\bar{\Omega}^t=\bar{\Omega}_1^t\cup\bar{\Omega}_2^t$, 
is defined as 
\[{\mathcal A}^t = {\mathcal A}(x, t) := x + d (x, t)\] 
for all $x\in \Omega_R$ and $t\in I$. 
In our model problem, 
we shall interprete $d$ as 
the finite element mesh movement, which defines 
the change of the computational sub-domains and is explicitly precomputed. 
The ALE time derivate of the function $u:\Omega^t\mapsto{\mathbb R}^3$ is defined 
as \[\partial_t u|_{{\mathcal A}^t}:=\partial_t u({\mathcal A}^t(x,t) , t)\] 
for all $x\in \Omega_R$ and $t\in I$. By the chain rule, we obtain 
\[\partial_t u=\partial_t u|_{{\mathcal A}^t}-w\cdot\nabla u\]
with $w=\partial_t {\mathcal A}^t\circ {{\mathcal A}^t}^{-1}$. 
Then we have the following model problem under the ALE framework: 
Find the solution $u : \Omega^t \mapsto {\mathbb R}^3$, for all 
$t\in I $, such that
\begin{equation}\label{eq:heatunfixedinterface}
  \begin{aligned}
    \partial_t u|_{{\mathcal A}^t} - w\cdot\nabla u - \nabla\cdot (a \nabla u ) = 0 & \textup{ in } \Omega_1^t\cup\Omega_2^t ,\\
            u_1 = u_2  & \textup{ on } \Gamma^t, \\
            a_1\frac{\partial u_1}{\partial n_1} + a_2\frac{\partial u_2}{\partial n_2}   
            & \textup{ on } \Gamma^t \\
  \end{aligned}
\end{equation}
with the initial conditions $u(x,0)=0$, $w(x,0)=0$ for all 
$x\in\Omega_1^0\cup\Omega_2^0$, and 
the boundary conditions $u = g_D$ on $\Gamma_D$ and 
$a\frac{\partial u}{\partial n}  = g_N$ on $\Gamma_N$ at $t>0$. 
Here 
$a=a_1\in {\mathbb R}^{+}$ in $\Omega_1^t$, 
$a=a_2\in {\mathbb R}^{+}$ in $\Omega_2^t$, $a_1\neq a_2$, are two different material 
coefficients in two moving domains, respectively.

\subsection{A combination of the ALE and macro-element method}
\begin{figure}[htbp]
  \centering
  \scalebox{0.7}{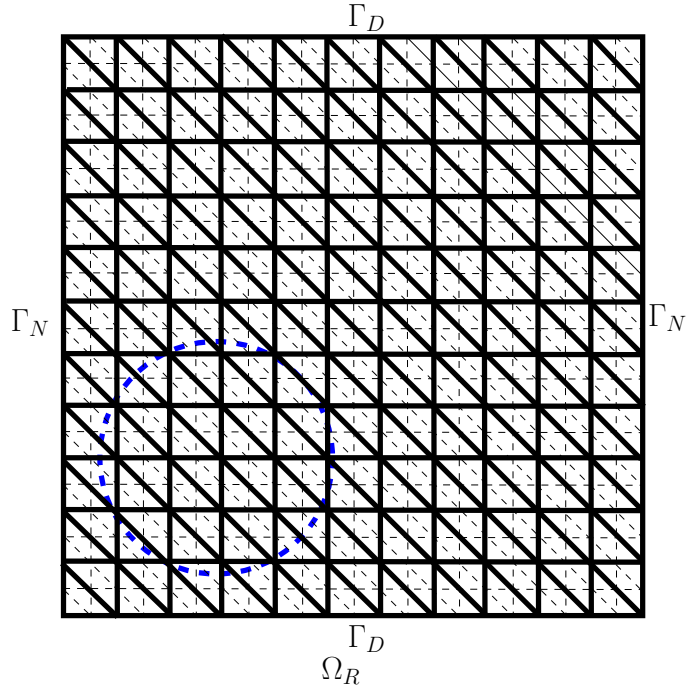}
  \caption{An illustration of a reference domain $\Omega_R$ decomposed into 
    macro-elements: macro-element edges (thick solid lines), 
    introduced new sub-element edges (thin dashed lines) and the interface (thick blue dashed line).}
  \label{fig:aleref}
\end{figure}

In the classical ALE method, we use the interface tracking method, 
where the mesh grids on the 
interface are following the object movement. The mesh movement 
inside the computational domains is computed by an arbitrary extension into the domain, 
e.g., a simple harmonic extension. The main drawback of this method is 
the restriction to small deformations. In case of large deformation or displacement, 
the mesh quality may deteriorate rapidly. To overcome this difficulty, we develop an 
interface capturing method, that is a combination of the ALE and 
macro-element method \cite{RBM11, YW11:00}. According to the cutting cases, 
the underlying reference domain is decomposed into macro-elements: 
four triangles in each macro-element in 2D and 
four tetrahedra plus one octahedron in 3D, see Fig.~\ref{fig:aleref} for an illustration 
of such decomposed reference domain into structured grids in 2D. 
The velocity $w:\Omega_R\mapsto {\mathbb R}^3$ of the mesh movement  
is constructed locally in each sub-element of the macro-element by an interpolation. 
The same applies to the displacement $d:\Omega_R\mapsto {\mathbb R}^3$ 
of the mesh movement, with respect to the reference configuration $\Omega_R$. 
The local velocity and displacement are related by $w=\partial_t d$. 
We comment that, for cells that are completely 
untouched with the moving interface (far away from the moving object), 
the velocity $w=0$ and the 
ALE mapping is an identity. 
In this case, the equation (\ref{eq:heatunfixedinterface}) is reduced to the one under the 
usual Eulerian framework. 

\begin{figure}[htbp]
  \centering
  \scalebox{1.2}{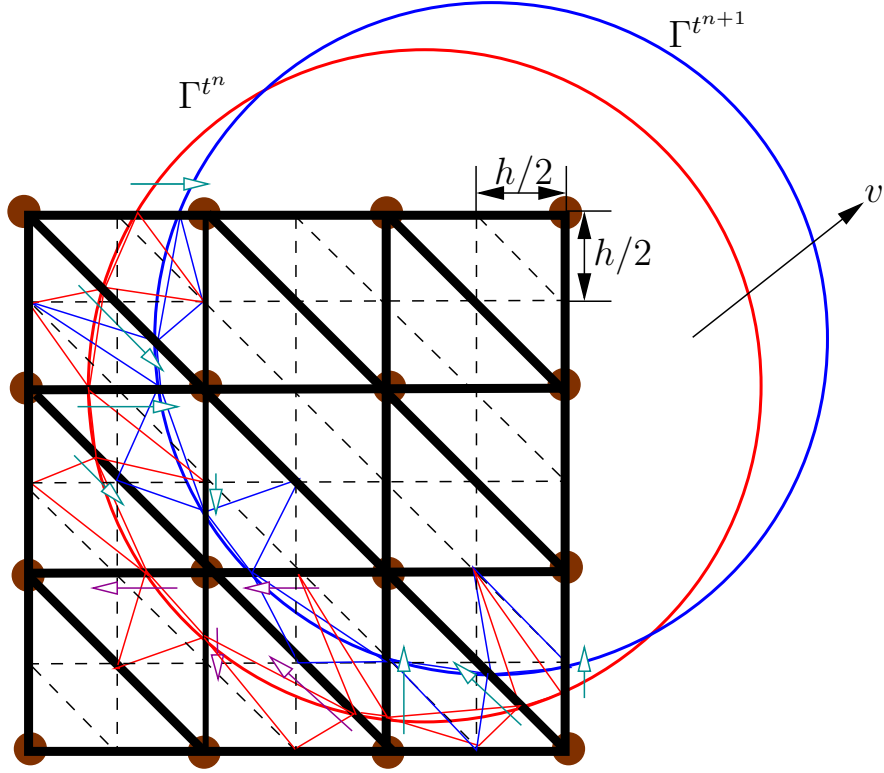}
  \caption{An illustration of the local nodes movement with restriction to each macro-element: 
  fixed macro-element nodes (brown dots), reconstructed moving 
  interface and locally adapted 
  triangle mesh at $t=t^n$ with nodes from the intersection (red lines), 
  reconstructed moving interface and locally adapted triangle mesh at $t=t^{n+1}$ 
  with nodes from the intersection (blue lines),  the moving direction of the intersection nodes within 
  each macro-element at the interface, cyan arrows (none of the intersection nodes is the edge middle point), 
  magenta arrows (one of the nodes is the edge middle point). }
  \label{fig:ale}
\end{figure}

\section{Temporal and spatial discretization}\label{sec:dis}
\subsection{Temporal discretization}
Let the time interval $I$ be divided into $N$ equidistant 
small time intervals $\Delta t$, i.e., $\Delta t = T/N$.  Let $t^n=n\Delta t$ be 
the time at level $n$. By the notation $f^n=f(x, t^n)$, 
we denote the function defined at the time $t^n$ and in the 
corresponding domain. We employ first-order implicit 
Euler scheme to discretize the time derivative: For all $n\geq 1$ and given 
$u^0=u_0$, we have 
\begin{equation}\label{eq:timedis}
  \begin{aligned}
  \frac{u^n - \tilde{u}^{n-1}}{\Delta t}
  -w^n\cdot\nabla u^n - \nabla\cdot(a\nabla u^n) = 0 & \textup{ in } \Omega_1^{t^n}
  \cup\Omega_2^{t^n} ,\\
     u_1^n = u_2^n  & \textup{ on } \Gamma^{t^n}, \\
            a_1\frac{\partial u_1^n}{\partial n_1^n} + a_2\frac{\partial u_2^n}{\partial n_2^n}   
            & \textup{ on } \Gamma^{t^n}, \\
  \end{aligned}
\end{equation}
where $\tilde{u}^{n-1}=u^{n-1}\circ ({\mathcal A}^{t^{n-1}})\circ ({{\mathcal A}^{t^n}})^{-1}$,
$w^n=\frac{d^n-d^{n-1}}{\Delta t}\circ ({\mathcal A}^{t^n})^{-1}$. In our model problem, 
the displacement $d$ is defined on the reference domain $\Omega_R$ and is explicitly 
evaluated by the intersection of the moving object with the underlying tetrahedral mesh. In 
Fig.~\ref{fig:ale}, the cyan and magenta arrows indicate the movement of the 
intersection points with respect to the underlying reference macro-element mesh. Assume 
the underlying mesh size consisting 
of the macro-element is $h$, then the mesh movement velocity is controlled by 
$|w|\approx \frac{h}{2\Delta t}$. When we choose sufficiently small mesh size, it will 
only introduce a very small convection term in the model problem, that is in general 
less problematic to perform standard finite element discretization and 
to solve the arising linear system of equations. 

\subsection{Spatial discretization}
The weak formulation arises from (\ref{eq:timedis}) by integration by parts and 
reads as follows: Find the solution 
$u^n\in V_g=\{g_D + V_0\}$ with $V_0=H_0^1(\Omega)^3=\{v\in H^1(\Omega)^3 | v=0 \text { on }\Gamma_D\}$ 
such that, for all $v\in V_0$, we have
\begin{equation}\label{eq:weakform}
  \left(\frac{u^n-\tilde{u}^{n-1}}{\Delta t}, v \right)_{\Omega^{t^n}}
  - \left(w^n\cdot\nabla u^n , v\right)_{\Omega^{t^n}}
  +(a\nabla u^n, \nabla v)_{\Omega^{t^n}} = \langle g_N, v\rangle_{\Gamma_N} ,
\end{equation}
where 
the continuity condition for the solution on $\Gamma^t$
has been explicitly enforced by using one identical 
$u^n$ in the domain $\Omega^{t^n}$, and the surface traction balance 
condition is implicitly included in the week form by integration by parts. 

\begin{figure}[htbp]
  \centering
  \scalebox{0.5}{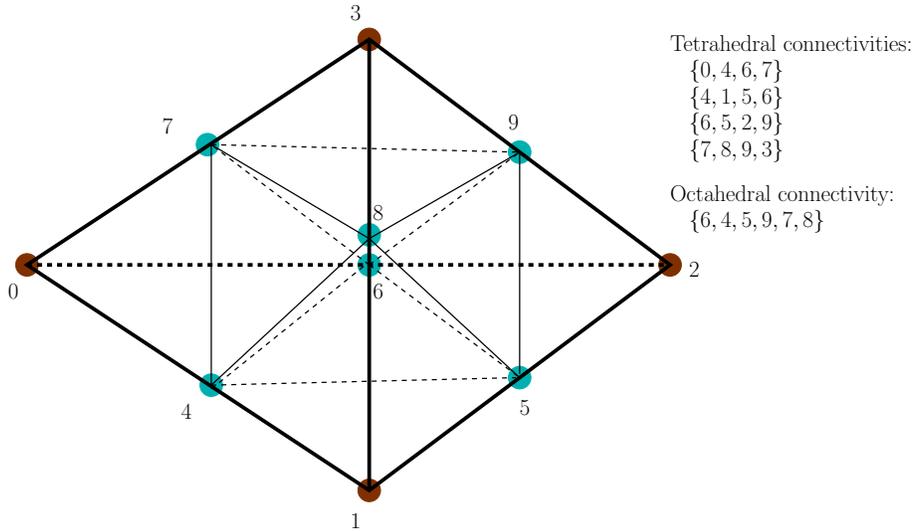}
  \caption{An illustration of a macro-element with five small sub-elements.}
  \label{fig:tetverts}
\end{figure}

\begin{figure}[htbp]
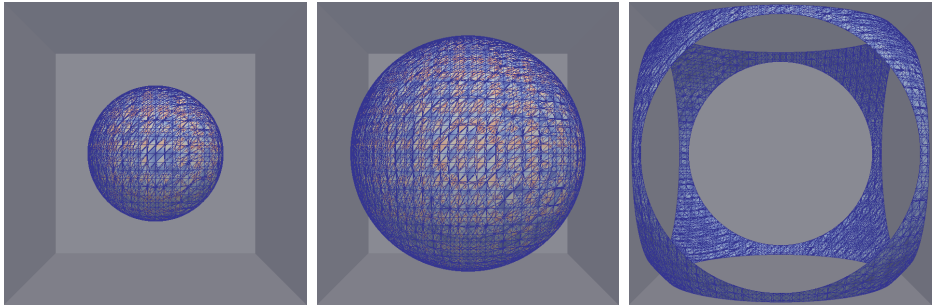

  \centering
  \includegraphics[scale=0.32]{constructface1.pdf} 
  \includegraphics[scale=0.32]{constructface2.pdf} 
  \includegraphics[scale=0.32]{constructface3.pdf} 
  \caption{A sequence of reconstructed surface meshes of the immersed growing objects.}\label{fig:surfconstruct}
\end{figure}

We use a finite element method for the spatial discretization. 
This method relies on the piecewise linear 
basis functions constructed on the underlying hybrid mesh consisting of 
tetrahedral and octahedral elements. Such mixed elements are obtained by 
decomposing each macro-element (a big tetrahedra) 
into four tetrahedral elements and one octahedral element; 
see Fig.~\ref{fig:tetverts} for an illustration of such a typical macro-element. 
Each tetrahedral macro-element has four fixed nodes with local node numbering $0-3$ 
(brown dots in Fig.~\ref{fig:tetverts}) 
and six nodes $4-9$ on edges (cyan dots in Fig.~\ref{fig:tetverts}) 
that are given by the edge middle points or the intersection points between the 
edge and the moving object. Each macro-element is decomposed into 
five sub-elements: 
four tetrahedron with the local node numbering $\{0, 4, 6, 7\}$, 
$\{4, 1, 5, 6\}$, $\{6, 5, 2, 9\}$, $\{7, 8, 9, 3\}$ 
and one octahedron with the local node numbering $\{6, 4, 5, 9, 7, 8\}$. 
This gives a very limited intersection patterns. In addition every macro-element has very similar structure to each other, that is easy to templatize on the 
computer implementation. By this means, we are able to reconstruct the triangle surface mesh 
of the immersed object; see Fig.~\ref{fig:surfconstruct} for an illustration 
of a sequence of such surface meshes. 
The hybrid mesh, consisting of different element types, has also been used recently 
in the cardiac electrophysiology simulation \cite{RBM11} and in the fluid-structure 
interaction simulation \cite{YW11:00,YH11:00}. 

\begin{figure}[htbp]
  \centering
  \scalebox{1.2}{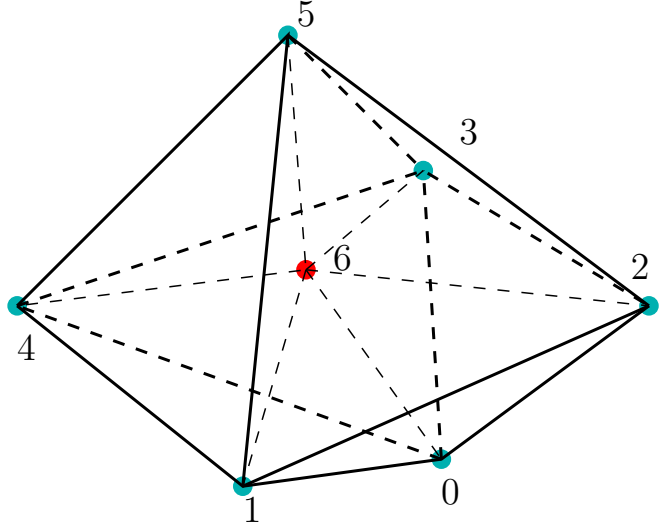}
  \caption{Splitting of an octahedron into $8$ tetrahedra $\{0,1,2,6\}$, 
    $\{0,2,3,6\}$, $\{0, 3, 4, 6\}$, $\{0, 1, 4, 6\}$, $\{5, 1, 4, 6\}$, $\{5, 1, 2, 6\}$, 
    $\{5, 2, 3, 6\}$, $\{5, 3, 4, 6\}$: Original edges (thick lines), added 
  edges (thin lines),  original nodes $\{0-5\}$, added node $\{6\}$.}
  \label{fig:octasplit}
\end{figure}

To be more precise, the finite element basis functions 
on the four tetrahedra in each macro-element is constructed as 
the standard hat function in 3D. On the remaining octahedron, 
we first add an auxiliary point $6$ near or at the mass center. The octahedron will be 
sub-divided into $8$ tetrahedra; see Fig.~\ref{fig:octasplit} for an illustration. 
We then construct standard hat functions on 
each tetrahedron. The extra degree of freedom at the node $6$ will be eliminated 
by the averaging of the values at nodes $0-5$; see more details in \cite{RBM11,YW11:00}. 
By this means, we do not introduce new degrees of freedom.  The number of 
total degrees of freedom is the number of nodes plus edges in the original mesh 
consisting of pure big tetrahedral macro-elements. 

\section{Solution methods for the linear system of equations}\label{sec:linsol}
\subsection{An all-at-once method}
After using finite element discretization, at each time step, we obtain the following 
linear system of equations:
\begin{equation}\label{eq:lineq}
  Ku=
  \left[\begin{array}{cc}
    A_{VV}& A_{VE}\\
    A_{EV} & A_{EE}
  \end{array}\right]
  \left[\begin{array}{c}
    u_V\\
    u_E
  \end{array}\right]=
  \left[\begin{array}{c}
    f_V\\
    f_E
  \end{array}\right]
  =f.
\end{equation}
We solve the linear system of equations by the AMG 
preconditioned conjugate gradient (PCG) 
method (see, e.g., \cite{UH:02}) and the AMG preconditioned 
GMRES method (see \cite{Saad86}). We mention here that, due to 
the small convection term, we found out that even the PCG method works well for solving 
such a non-symmetrically perturbed symmetric 
linear system of equations. For convenience of the solution 
procedure, the linear system has been ordered with firstly the degrees of freedom on the 
original tetrahedral nodes $u_V$, and then of the edges $u_E$, where the subscripts 
$V$ and $E$ are associated with the nodes and edges. Such reordering has been used in 
the AMG method for high-order finite element discretized equations \cite{YH14,ULHY15}. 
The stiffness matrices $A_{VV}$ and  $A_{EE}$ arise 
from the finite element assembly of the basis functions associated with the original macro-element 
nodes and edges, respectively, $A_{EV}$ and $A_{VE}$ are coming from the coupling. To solve 
such a linear system of equations, we use a special AMG method \cite{FK98:00}, that is 
based on the matrix graph connectivity. Similar idea was also developed in \cite{DB95}. 
In our numerical simulation, such solution methods give us quite 
satisfactory
results. 
We observe a quite robust behavior of the AMG preconditioner with respect to 
moving interface in each time step. 

\subsection{A segregated method}
\begin{figure}[htbp]
  \centering
  \includegraphics[scale=0.8]{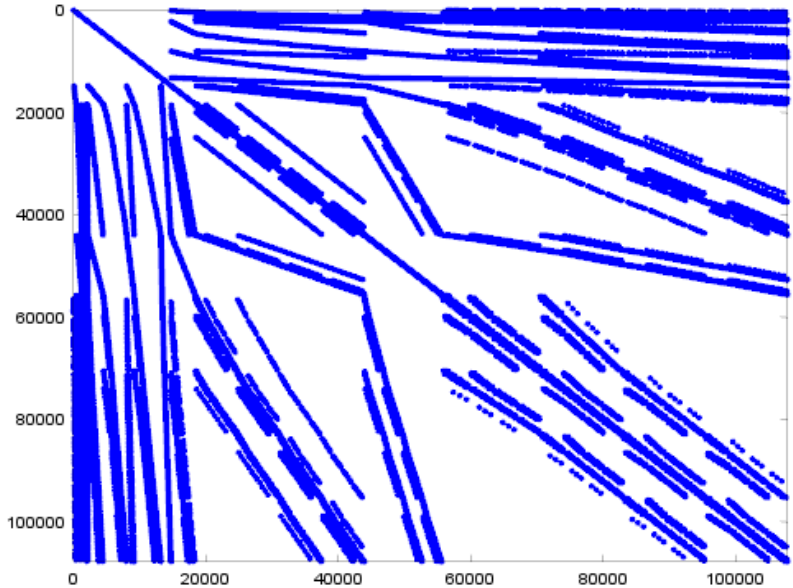} 
  \caption{Sparsity pattern of the system matrix $K$.}\label{fig:matrixpattern}
\end{figure}

By a close look at the matrix structure in (\ref{eq:lineq}), we have observed that 
$A_{VV}$ is a block-diagonal matrix. This is due to the fact that the degrees of 
freedom associated with the original macro-element nodes are 
completely decoupled. In Fig.~\ref{fig:matrixpattern}, we demonstrate a sparsity 
pattern of the system matrix $K$, where it is easy to see the block-diagonal 
structure of $A_{VV}$.

We now perform a $LU$ factorization of  the system matrix $K$ in (\ref{eq:lineq}):
\begin{equation}\label{eq:ldu}
  K=LU=
  \left[
    \begin{array}{cc}
      A_{VV} & 0\\ A_{EV} & S
    \end{array}\right]
  \left[
    \begin{array}{cc}
      I & A_{VV}^{-1}A_{VE} \\ 0 & I 
    \end{array}\right],
\end{equation}
where $S$ denotes the Schur complement $S=A_{EE}-A_{EV}A_{VV}^{-1}A_{VE}$. Since 
$A_{VV}^{-1}$ can be constructed very easily, 
the Schur complement $S$ can also be constructed exactly. A simple blockwise 
forward and backward substitution gives rise to the solution of the linear system. The 
main cost is to solve the Schur complement equation 
\begin{equation}\label{eq:schureq}
Sx_{EE}=b_{E} 
\end{equation}
for $b_{E}:=f_E-A_{EV}A_{VV}^{-1}f_V$. This is realized by applying the AMG preconditioned CG method \cite{UH:02}. 

\section{Numerical results}\label{sec:numa}
\subsection{The numerical result for the model problem 
  with an immersed moving sphere}
\begin{figure}[htbp]
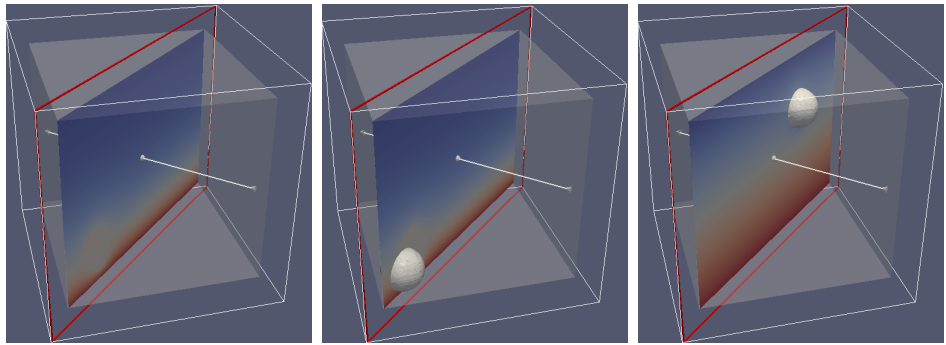

  \centering
  \includegraphics[scale=0.4]{statheat_cuttingface.pdf} 
  \includegraphics[scale=0.4]{statheat_cuttingface_t0.pdf}
  \includegraphics[scale=0.4]{statheat_cuttingface_t1.pdf}
  \caption{Cutting plane (left), constructed moving sphere surfaces at the time $t=0$ (middle) 
    and $t=0.5625$ (right).}\label{fig:heatmovingobjcutting}
\end{figure}

\begin{figure}[htbp]
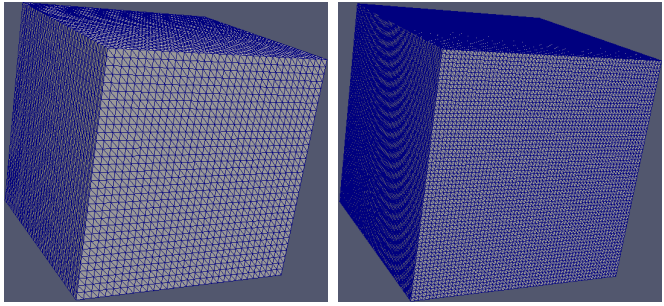

  \centering
  \includegraphics[scale=0.4]{meshorg.pdf} 
  \includegraphics[scale=0.4]{meshcut.pdf}
  \caption{Original pure tetrahedral mesh (left), the sub-divided mesh (right).}\label{fig:meshinfo}
\end{figure}

In the first example, we consider a sphere with fixed 
radius $0.12$ and the initial 
center at $(0.125, 0.125, 0.125)$ immersed in a unit cube; see 
Fig. \ref{fig:heatmovingobjcutting} for an illustration.  
The cube is decomposed into macro-elements with $35937$ nodes and 
$196608$ tetrahedra. The sub-divided hybrid mesh consists of $274625$ nodes and 
$786432$ tetrahedra and $196608$ octahedra. The total number of degrees of freedom is 
$823875$. See Fig. \ref{fig:meshinfo} for an illustration. 

\begin{figure}[htbp]
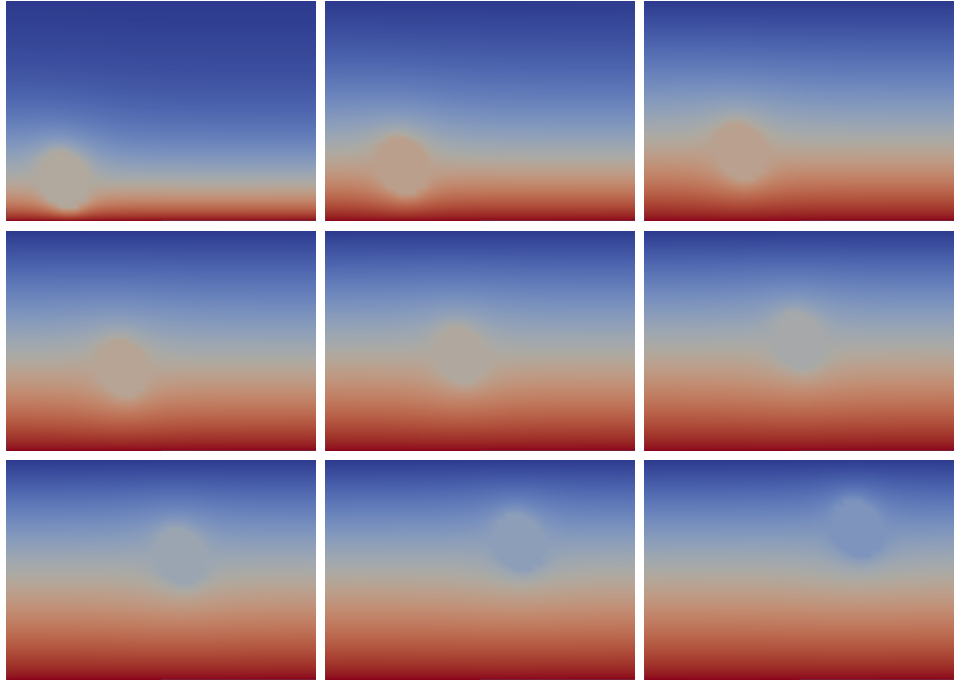

  \centering
  \includegraphics[scale=0.23]{heat1.pdf}
  \includegraphics[scale=0.23]{heat2.pdf}
  \includegraphics[scale=0.23]{heat3.pdf}\\[0.1cm]
  \includegraphics[scale=0.23]{heat4.pdf}
  \includegraphics[scale=0.23]{heat5.pdf}
  \includegraphics[scale=0.23]{heat6.pdf}\\[0.1cm]
  \includegraphics[scale=0.23]{heat7.pdf}
  \includegraphics[scale=0.23]{heat8.pdf}
  \includegraphics[scale=0.23]{heat9.pdf}
  \caption{Simulation results with the moving object in the domain 
    at different time levels $t=0.0625k$, $k=1,...,9$, on the cutting face.}\label{fig:heatmovingobj}
\end{figure}

\begin{figure}[htbp]
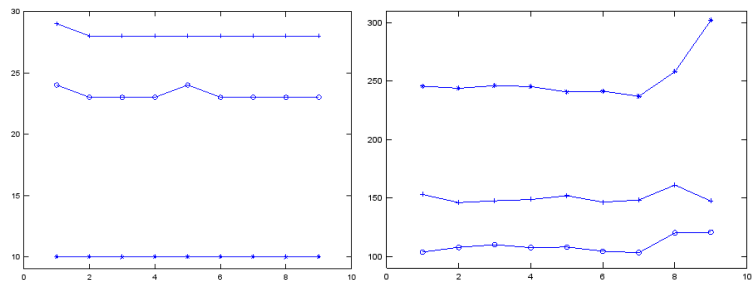

  \centering
  \includegraphics[scale=0.5]{heatamg.pdf}
  \includegraphics[scale=0.55]{heattimeamg.pdf}
  \caption{Iteration numbers (left) and CPU time measured in second $s$ (right) for solving 
    the time dependent heat equation with the immersed moving object in each time step: 
    AMG preconditioned CG (solid lines with circle markers), AMG preconditioned GMRES 
    (solid lines with star markers) in the monolithic method, 
    AMG preconditioned CG (solid lines with plus markers) 
    for the Schur complement equation in the segregated method.}\label{fig:heatamg}
\end{figure}

The sphere is moving along 
the line with the starting point $(0,0,0)$ and the ending point $(1,1,1)$, 
and the moving speed is $v=(1,1,1)^T$. The constructed sphere surface 
is shown in the middle and right plots of Fig.~\ref{fig:heatmovingobjcutting}. 
On the bottom of the cube, we set the Dirichlet boundary condition 
$u=(0, 0, 0)^T$, on the top, $u=(1, 0, 0)^T$. For the rest of the boundaries, 
we use the homogeneous Neumann boundary condition. The time stepsize is 
$\Delta t =0.0625$ and the number of time steps is $9$, i.e., the ending time is 
$T=0.5625$. The material coefficient inside $\Omega_1^t$ is $a_1=1.0e+06$ 
and inside $\Omega_2^t$ is $a_1=1.0$. The simulation results 
at different time on the cutting plane (see the left plot in Fig.~\ref{fig:heatmovingobjcutting}) 
are shown in Fig.~\ref{fig:heatmovingobj}. The relative residual error is set to $10e-09$ as stopping 
criteria. The iteration numbers and the computational CPU time (in second) 
of the AMG preconditioned CG and GMRES methods in the all-at-once method, and 
the iteration numbers of the AMG preconditioned CG for the Schur complement equation and 
the computational CPU time in the segregated method, 
are shown in Fig.~\ref{fig:heatamg}. We observe that, in terms of iteration numbers, the 
GMRES method shows the best performance, then the CG method, and last the segregated method. 
However, regarding CPU time, we see that, the CG method shows its best performance, then 
the segregated method, and last the GMRES. 

\subsection{The numerical result for the model problem 
  with an immersed growing sphere}
\begin{figure}[htbp]
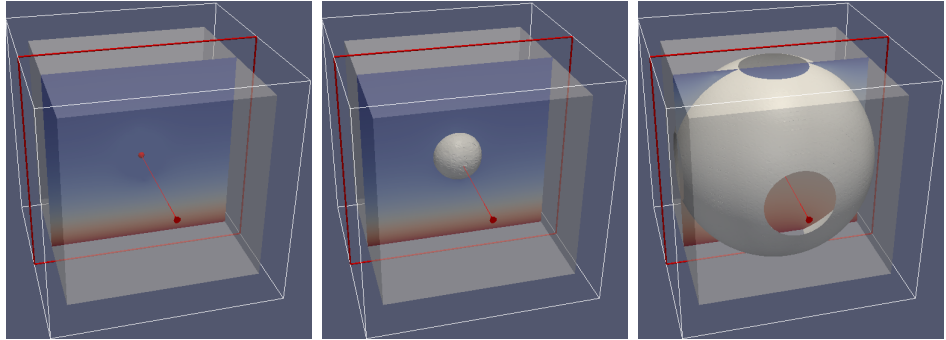

  \centering
  \includegraphics[scale=0.4]{statheatgrow_cuttingface.pdf} 
  \includegraphics[scale=0.4]{statheatgrow_cuttingface_t4.pdf}
  \includegraphics[scale=0.4]{statheatgrow_cuttingface_t8.pdf}
  \caption{Cutting plane (left), constructed growing sphere surfaces at time $t=0$ (middle) 
    and $t=0.45$ (right).}\label{fig:heatgrowingobjcutting}
\end{figure}

In the second example, we consider a sphere with an initial 
radius $0.08$ and the initial 
center at $(0.5, 0.5, 0.5)$ immersed in a unit cube; see 
Fig.~\ref{fig:heatgrowingobjcutting} for an illustration. We use the same finite element mesh 
as in the first example. 
The sphere is growing along the radius 
direction and the growing speed is $v=n$, where $n$ denotes 
the outward unit normal vector in the radius direction. The surfaces of the growing 
sphere at time $t=0$ and $t=0.45$ are constructed as shown in the middle and right plots of 
Fig.~\ref{fig:heatgrowingobjcutting}, respectively. 
On the bottom of the cube, we set the Dirichlet boundary conditions 
$u=(0, 0, 0)^T$, on the top, $u=(1, 0, 0)^T$. For the rest of the boundaries, 
we use the homogeneous Neumann boundary condition. The time stepsize is 
$\Delta t =0.05$ and the number of time steps is $9$, i.e., the ending time is 
$T=0.45$. The material coefficient inside $\Omega_1^t$ is $a_1=1.0e+06$ 
and inside $\Omega_2^t$ is $a_1=1.0$. The simulation results on the cutting plane (see the left plot in 
Fig. \ref{fig:heatgrowingobjcutting}) is shown in Fig.~\ref{fig:heatgrowingobj}. 
The relative residual error is set to $10e-09$ as stopping criteria of the linear solvers. 
The iteration numbers and the computational CPU time (in second) 
of the AMG preconditioned CG and GMRES methods in the all-at-once method, and 
the iteration numbers of the AMG preconditioned CG for the Schur complement equation and 
the computational CPU time in the segregated method, 
are shown in Fig.~\ref{fig:heatgamg}. We observe that, 
in terms of iteration numbers, the 
GMRES method shows the best performance, then the CG method, and last the segregated method. 
However, regarding CPU time, we see that again, the CG method shows the best performance, then 
the segregated method, and last the GMRES method. 
\begin{figure}[htbp]
  \centering
  \includegraphics[scale=0.32]{heatg1.pdf}
  \includegraphics[scale=0.32]{heatg2.pdf}
  \includegraphics[scale=0.32]{heatg3.pdf}\\[0.1cm]
  \includegraphics[scale=0.32]{heatg4.pdf}
  \includegraphics[scale=0.32]{heatg5.pdf}
  \includegraphics[scale=0.32]{heatg6.pdf}\\[0.1cm]
  \includegraphics[scale=0.32]{heatg7.pdf}
  \includegraphics[scale=0.32]{heatg8.pdf}
  \includegraphics[scale=0.32]{heatg9.pdf}
  \caption{Simulation results with the growing object in the domain 
    at different time levels $t=0.05k$, $k=1,...,9$, on the cutting face.}\label{fig:heatgrowingobj}
\end{figure}
\begin{figure}[htbp]
  \centering
  \includegraphics[scale=0.5]{heatgamg.pdf}\;
  \includegraphics[scale=0.55]{heatgtimeamg.pdf}
  \caption{Iteration numbers (left) and CPU time measured in second $s$ (right) for solving 
    the time dependent heat equation with the immersed growing object in each time step: 
    AMG preconditioned CG (solid lines with circle markers), AMG preconditioned GMRES 
    (solid lines with star markers) in the monolithic method, 
    AMG preconditioned CG (solid lines with plus markers) 
    for the Schur complement equation in the segregated method.}\label{fig:heatgamg}
\end{figure} 

\section{Conclusion}\label{sec:con}
In this work, we develop an ALE method on the underlying 
reference domain decomposed macro-elements consisting of 
tetrahedral and octahedral elements. 
That is combined with the interface capturing method. 
The numerical results demonstrate the robustness 
of this method with respect to large displacement or deformation of the 
moving interface in the model parabolic problem. We have compared the algebraic multigrid based 
all-at-once and the segregated methods for solving the linear system of algebraic equations arising from 
the finite element discretization. We observed that the all-at-once AMG preconditioned CG method shows 
the best performance in terms of CPU time. The segregated method shows comparable performance. 
Regarding the iteration numbers, the AMG preconditioned GMRES method shows the best performance. 


\bibliography{interface}
\bibliographystyle{siam}

\end{document}